\documentclass[10pt,twocolumn]{article}
\usepackage{setspace}

\usepackage[english]{babel}

\usepackage{amsmath,amssymb,amsthm,amsfonts}
\usepackage{graphicx}
\usepackage{booktabs}

\usepackage{chngcntr}
\usepackage{textcomp}
\usepackage[hyphens]{url}
\usepackage{enumitem}
\usepackage{rotating}
\usepackage{caption}
\usepackage{float}
\usepackage{mathrsfs}
\usepackage{mathtools}
\usepackage{tabularx}
\usepackage{bm}
\usepackage[export]{adjustbox}
\usepackage{subcaption}
\usepackage{siunitx}
\usepackage{wrapfig}
\usepackage{hhline}
\usepackage{array}

\usepackage[authoryear,round]{natbib}

\usepackage[table]{xcolor}
\usepackage{soul}
\definecolor{dr}{rgb}{0.75,0.00,0.00}
\definecolor{lr}{rgb}{1.00,0.75,0.75}
\sethlcolor{lr}

\usepackage[colorinlistoftodos]{todonotes}

\newcommand{\aiopt}{\texttt{AI4OPT}}

\title{\aiopt{}: AI Institute for Advances in Optimization}
\author{Pascal Van Hentenryck and Kevin Dalmeijer \\ Georgia Institute of Technology \\ \{pvh, dalmeijer\}@gatech.edu}
\date{}

\graphicspath{{./figures/}}

\usepackage[colorlinks=true, allcolors=blue]{hyperref}
\usepackage[letterpaper, top=0.97in, bottom=0.97in, left=1.0in, right=1.0in]{geometry}

\begin{document}

\twocolumn[\vspace{-2em plus 1em minus 0.5em}\maketitle]

\begin{figure}[!t]
	\centering
	\includegraphics[width=0.9\linewidth,trim=0.5cm 0 0 0,clip]{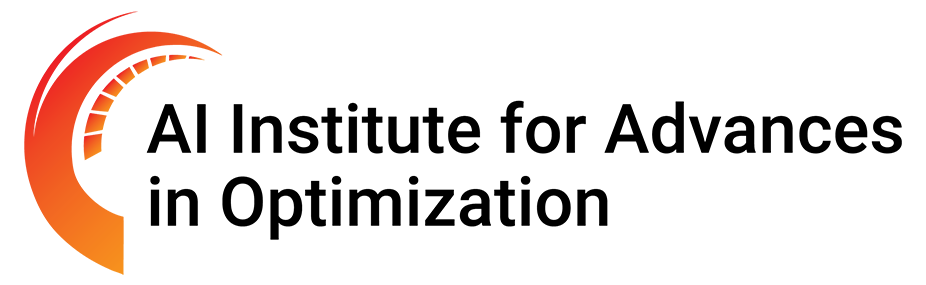}
\end{figure}

\begin{abstract}
  This article is a short introduction to \aiopt{}, the NSF AI Institute for Advances in Optimization.
  \aiopt{} fuses AI and Optimization, inspired by end-use cases in supply chains, energy systems, chip design and manufacturing, and sustainable food systems. \aiopt{} also applies its ``teaching the teachers'' philosophy to provide longitudinal educational pathways in AI for engineering.
	\begin{flushleft}\emph{\textbf{Keywords}:
		optimization, machine learning,\\supply chains, energy systems, chip design and manufacturing, resilience, sustainability.
		% up to 7 keywords
	}\end{flushleft}
\end{abstract}

\section{Introduction}

The mission of the NSF Artificial Intelligence (AI) Institute for
Advances in Optimization (\aiopt{}) is
\begin{quote}
{\em to revolutionize decision making at massive scales by fusing AI and mathematical optimization and delivering scientific breakthroughs that the two fields cannot achieve
 independently.}
\end{quote}
This Institute pursues this objective by integrating the model-driven
paradigm typically followed in Operations Research with data-driven
methodologies coming from AI. The research in \aiopt{} is
use-inspired, addressing fundamental societal and technological
challenges of our times. They include 
\begin{itemize}[noitemsep,topsep=2.5pt,parsep=0pt,partopsep=0pt,leftmargin=0.5cm]
\item How to design agile, sustainable, resilient, and equitable supply chains?	
\item How to operate energy systems powered by distributed renewable energy resources?
\item How to deliver a step change in chip design and manufacturing, and manufacturing as a whole?
\item How to create sustainable eco-systems within the food-water-energy
  nexus?
\end{itemize}
\aiopt{} focuses on {\em AI for Engineering}, which raises deep
scientific challenges in terms of reliability, robustness, and
scalability. Indeed, \aiopt{} is driven by high-stake applications
that feature physical, engineering, and business constraints, and
complex objectives that must balance efficiency, resilience,
sustainability, and equity. Moreover, the underlying optimization
problems at the core of the grand challenges are of very large scale,
many of which are beyond the scope of existing technologies. To
address these, \aiopt{} is organized around methodology thrusts that
focus on specific challenges: they include a new generation of
data-driven optimization solvers, decision making under uncertainty,
combinatorial and reinforcement learning, end-to-end optimization, and
decentralized learning and optimization. In addition, a transversal
thrust on Ethical AI ensures that ethics is included in the design of
every fundamental and use-inspired project, not as an afterthought.
The complementarity of the end-use cases and the methodology thrusts
creates a virtuous cycle of innovation, both in foundational research
and industrial impact.

The research mission of \aiopt{} is complemented by its educational vision which is
\begin{quote}
  {\em to create longitudinal pathways for AI in engineering, from high-school to graduate education, using a ``teach the teachers'' philosophy to maximize impact.}
\end{quote}
The pathways start in middle and high-schools (through summer camps and engineering practices), move to undergraduate education
through the Faculty Training Program, and the creation of graduate
programs. The {\em ``teaching the teachers''} philosophy is pervasive
across the Institute: its goal is to empower teachers at every
education level to create programs in their own institutions, e.g.,
minors and majors in artificial intelligence. \aiopt{}
also includes mentorship programs to help students take leadership
roles and reinforce their understanding of the material over time.
The Institute has a strong focus on historically black high schools,
Historically Black Colleges and Universities (HBCUs), and Minority
Serving Institutions (MSI), to develop talent and increase the
diversity of the AI workforce.

\aiopt{} is led by the Georgia Insitute of Technology in collaboration
with universities in California (UC Berkeley, USC, UC San Diego),
Texas (UT Arlington), and Georgia (Clark Atlanta University). The
Institute is creating a vibrant nexus in AI and optimization, bringing
together academic institutions, industrial partners, international
collaborators, and educators. In Atlanta, \aiopt{} is located on the
12th floor of CODA building in Midtown, providing a prime space for
faculty, research scientists, and students that encourages knowledge
cross-fertilization.

The Industrial Partner Program (IPP) of AI4OPT features novel
internship programs to facilitate research collaborations between
academia and industry. It assembles some of the most innovative
companies in supply chains, manufacturing, and energy systems, with
the goal of maximizing the impact emerging from the fusion of AI and
optimization.

The rest of this article is organized as follows. Section
\ref{section:proxies} illustrates how the Institute contributes
scientific breakthroughs that the two fields cannot achieve
independently. Section \ref{section:challenges} describes some of the
societal and technological challenges that are driving \aiopt{}
research. Section \ref{section:methodology} briefly reviews the
methodology thrusts driven by the end-use cases. Sections
\ref{section:education} and \ref{section:nexus} outline some workforce
development activities and how \aiopt{} acts as a nexus at the
intersection of AI and optimization. It is impossible to do justice
to all the activities of the Institute in a short article, but the
hope is that this presentation will encourage readers to learn more
about \aiopt{}.

\vspace{-0.1cm}
\section{Optimization Proxies}
\label{section:proxies}

At its core, optimization models are used in decision-making applications to map problem inputs into optimal solutions. Optimization solvers have seen dramatic
progress over the last decades, producing optimal solutions to many
industrial problems. For instance, optimization models quite literally keep the lights on by committing generators and dispatching electricity in real time
every five minutes. Optimization models run end-to-end supply chains, which involve aspects such as the scheduling of manufacturing plants, load consolidation  for middle-mile logistics, 
and the design of e-commerce
networks, to name a few. Yet recent developments are
challenging even the best solvers:
optimization models have grown even larger, are expected to capture the realities of an increasingly uncertain and volatile world, and are blurring the distinction between planning and operations. As a result, 
optimization solvers have become too slow in many
contexts: they include real-time applications, large-scale Monte-Carlo
simulations that are based on optimization, and environments where
humans interact with optimization technology.

In those circumstances, it is natural to explore whether machine
learning can replace optimization, moving most of the computational
burden offline. A machine-learning model can then learn the
input-output mapping of the optimization model, producing a first
approximation to the concept of an {\em optimization proxy}.
The challenge,
however, comes from applying this idea to {\em AI for
  Engineering}. Indeed, many of the end-use cases of the Institute
feature optimization problems with hard physical, engineering, and
business constraints. For instance, in an electrical grid, the load
(demand) and the generation (supply) must be equal at all times. In
supply chains, shipments must fit within the vehicle capacity.
In addition, optimization proxies will be deployed in
high-stake applications, which means that they must cater to a wide variety of instances, deliver high-quality solutions with performance guarantees, and learn
instances with millions of input parameters and make hundreds of thousands of predictions.

\begin{figure}[!t]
	\centering
	\includegraphics[width=\linewidth]{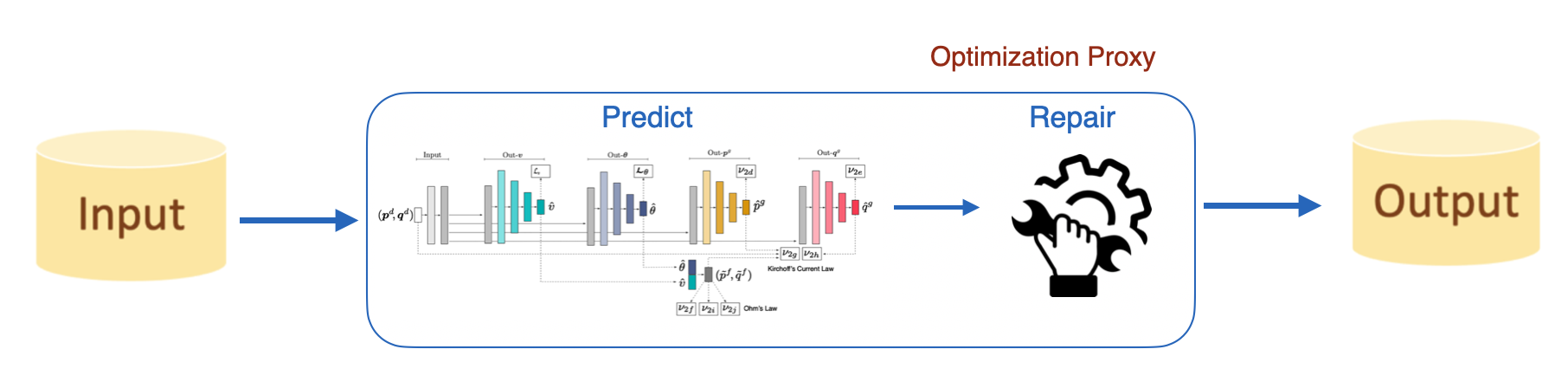}
	\caption{The Architecture of Optimization Proxies.}
        \label{fig:proxies}
\end{figure}

\begin{figure}[!t]
	\centering
	\includegraphics[width=\linewidth]{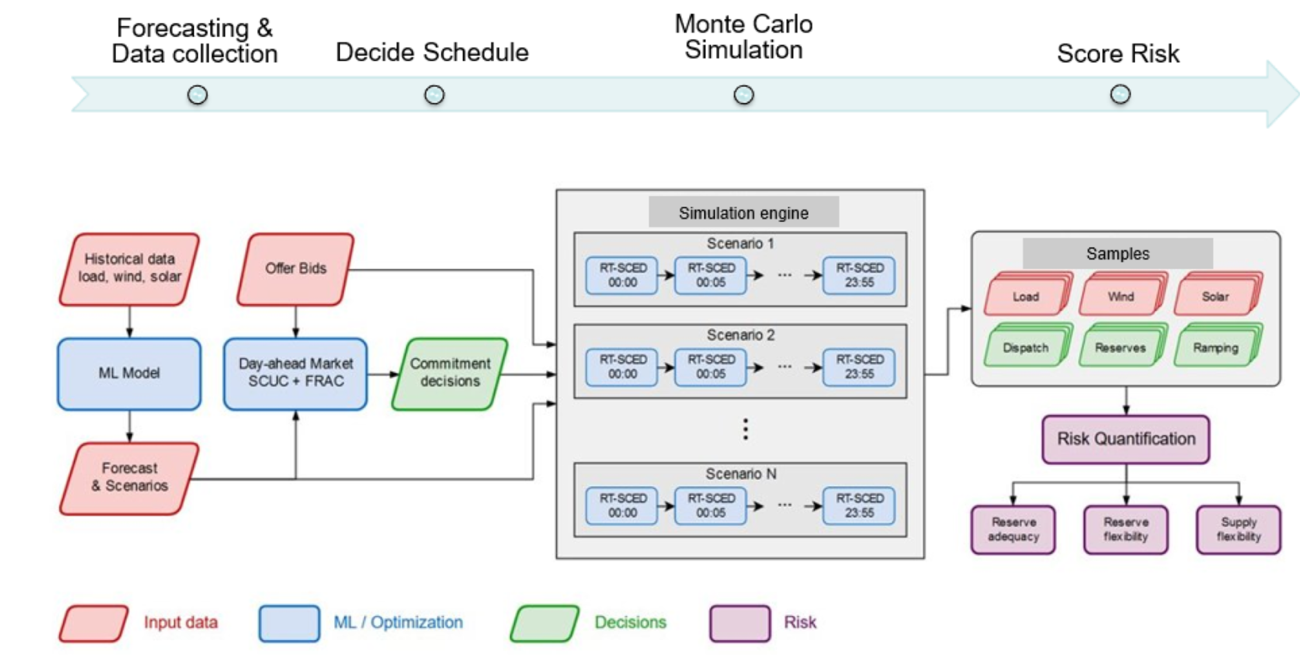}
	\caption{Real-Time Risk Assessment.}
        \label{fig:risk}
\end{figure}

To address these considerations, the {\em End-to-End Optimization}
thrust explores the science and engineering of optimization
proxies. Its research led to novel architectures including the one
depicted in Figure \ref{fig:proxies}. This architecture postulates an
optimization proxy as the composition of a machine-learning layer that
produces a high-quality approximation of the optimization model,
followed by a feasibility layer that repairs the prediction to deliver
a feasible solution. This Learning and Repair architecture can be
trained end-to-end, backpropagating the loss function through the
feasibility layers.

This fusion of AI and optimization can deliver breakthroughs that cannot
be achieved by the two fields independently. One such application is
real-time risk assessment, which is becoming increasingly pervasive in
energy systems and supply chains, including for some of the key
partners of \aiopt{}. As mentioned earlier, transmission
system operators typically solve a market-clearing optimization every
five minutes to balance load and generation.  Given the increased
volatility in net load due to intermittent renewable resources, grid
operators are interested in real-time risk assessment tools,
like the one described in Figure \ref{fig:risk}.  Such tools run a
large number of Monte-Carlo simulations, using scenarios from
probabilistic forecasts, to quantify the system-wide risk.  A single
simulation may take up to 45 minutes, given the computational
complexity and the sheer number of optimization problems, making it
impractical to assess risk in real time.  Work by the Institute has
shown that optimization proxies are a transformative technology for
this application: they make it possible to run these simulations in a
few milliseconds with relative errors below 1\%, giving rise to
potential new tools to manage risk in real time
\citep{ChenEtAl2022-LearningOptimizationProxies}.

\vspace{-0.1cm}

\section{End-Use Cases}
\label{section:challenges}

As mentioned, already, at \aiopt{}, challenges from end-use cases inspire foundational research, which then delivers innovations to address them. Here is a brief description of these end-use cases.

\paragraph{Supply Chains}

Supply-chain management used to be an arcane topic, discussed by a few
and invisible to the general public. This has changed after the
pandemic: the public is now aware of a topic that has become
top-of-mind in many corporate boards. Supply chains have become
larger, and e-commerce has proliferated, imposing significant
environmental costs to meet new customer expectations. At the same time,
many customers and suppliers, especially in rural regions, face
increasing difficulties in procuring or delivering specific
products. {\em What is needed is a paradigm change, a new vision for
  supply chains that complements efficiency with resilience,
  sustainability, and equity goals.}  Research in supply chains at
\aiopt{} is centered around end-to-end supply chains, with
scalability, resilience, sustainability, and equity as core
challenges. \aiopt{} has assembled a consortium of partners that cover
(almost) all aspects of supply chains. It leverages novel forecasting
methods, optimization proxies, decision making under uncertainty, and
automation to meet these challenges.

\paragraph{Energy Systems}

The challenge for energy systems is clear: {\em how to reinvent the planning and operations of a grid powered by renewable energy sources and
  storage}. Energy systems are transitioning from the century-old
``generation follows the load'' organization of the grid to a paradigm
centered on {\em risk assessment and risk management}. This end-use case 
helps reinvent energy systems by pursuing four overarching
themes: (1) probabilistic forecasting to quantify uncertainty; (2)
stochastic and risk-aware optimization to capture this uncertainty in
decision processes; (3) optimization proxies to perform real-time risk
assessment and risk-aware optimization; and (4) decentralized
optimization to address the massive proliferation of distributed
energy resources.

\paragraph{Chip Design and Manufacturing}

Each generation of chips is becoming more expensive to design,
requiring numerous cycles between expert designers and simulators. It
is no longer possible or desirable to separate the various phases of
the design, e.g., circuit synthesis, placement, and routing. {\em What is
needed is a new generation of tools that help engineers design
circuits more holistically.} Machine learning (including RL and inverse
learning) has been shown to have significant potential in this area; it is the goal of this end-use case to explore its role in circuit
optimization. Chip manufacturing has also becoming increasingly
complex and subject to uncertainty in supply, demand, and the
complexity of the bill of materials.
In conjunction with supply chains, this end-use case also aims at transforming the optimization of the chip manufacturing process.

\paragraph{Sustainable Systems}

The food-energy-water nexus is identified as one of the key grand
challenges of the 21st century and AI has demonstrated early potential
to address complex problems in this space. This end-use case conducts
three interconnected projects on biogas, water, and food to reduce greenhouse emissions and boost food production.

  \vspace{-0.1cm}
\section{Methodology Thrusts}
\label{section:methodology}

The methodology thrusts carry out foundational AI research inspired by the end-use cases. Here is a brief description of the methodology thrusts of \aiopt{}.

\paragraph{End-to-End Optimization} 

The End-to-End Optimization thrust primarily focuses on the {\em
  science and engineering of optimization proxies} that were described
in Section~\ref{section:proxies}. Recent contributions include the
concepts of self-supervised primal-dual learning
\citep{ParkVanHentenryck2023-SelfSupervisedPrimal}, compact learning
\citep{ParkEtAl2023-CompactOptimizationLearning}, and End-To-End
Learning and Repair \citep{ChenEtAl2023-EndEndFeasible}.  The thrust
draws inspiration from the end-use cases in energy systems and supply
chains, and also explores topics in decision-focused learning, learning to optimize, verification, explanation, and
formal guarantees.

\paragraph{New Generation Solvers}

The Solvers thrust works on a new generation of highly
tunable optimization solvers that use machine learning and historical
data to dramatically improve performance in settings where an
optimization model is used repeatedly. Recent results apply the
\textit{Learning to Optimize} paradigm to mixed-integer programming
\citep{HuangEtAl2023-SearchingLargeNeighborhoods}, mixed-integer
nonlinear programming \citep{FerberEtAl2023-SurcoLearningLinear}, and
AI planners \citep{HuangEtAl2023-DeadlineAwareMulti}.
The end-use cases contribute problem instances to benchmark solver performance.

\paragraph{Decision Making Under Uncertainty}

The energy and supply-chain end-use cases clearly indicate the need
for advances in decision making under uncertainty.
This includes probabilistic forecasting, uncertainty
quantification, scenario generation, and detection of rare events in
presence of spatial-temporal correlations
\citep{XuXie2022-SequentialPredictiveConformal}. Of particular
interest is the fundamental and applied research on conformal
predictions. The thrust also explores solution techniques for specific
classes of multi-stage stochastic optimization problems
\citep{LanShapiro2023-NumericalMethodsConvex} and new Bayesian
risk-sensitive and distributionally-robust optimization models
\citep{JuLan2023-DualDynamicProgramming}.

\paragraph{Reinforcement Learning}

The RL thrust focuses on the end-use cases of the Institute, which are
much larger and more complex than environments in which RL has been
successful so far.  It contributes foundational advances to deep RL to
handle such complex environments
\citep{LaskinEtAl2022-UnsupervisedReinforcementLearning}, and expands
RL research to include societal and ethical considerations.
The thrust will be increasingly focused on \emph{offline RL}
\citep{ChenMaguluri2022-SampleComplexityPolicy} to make the technology
safer and more amenable to industrial use.

\paragraph{Combinatorial Learning}

This thrust studies machine learning in the context of combinatorial
and highly constrained applications with the goal to improve
generalization and interpretability and reduce errors. Recent results
include highly-specific strong convex models with structured sparsity
\citep{AtamtuerkGomez2022-SupermodularityValidInequalities},
pairwise-based optimization algorithms for counteracting learning bias
\citep{HochbaumEtAl2023-BreakpointsBasedMethod}, and meta-algorithms
to automatically select the best solver
\citep{Asin-AchaEtAl2022-FastAlgorithmsCapacitated}.

\paragraph{Distributed and Multi-Agent Learning and Optimization}

This thrust explores decentralized solutions to manage a large number
of agents, motivated by applications in energy systems and automated
warehouses. Recent results focus on AI-based decentralized path
planning and execution \citep{XuEtAl2022-MultiGoalMulti}, on how
agents can help each other learn through communication
\citep{ZhangEtAl2022-DualAcceleratedMethod}, and
distributed learning algorithms that are robust against communication
imperfections \citep{ZengEtAl2023-ConnectedSuperlevelSet}.

\paragraph{Ethical AI}

The Ethical AI thrust is transversal: it draws from, and informs,
every thrust and end-use case in the Institute. It leverages the
fusion of optimization and machine learning to ensure {\it ethical and
  socially conscious design} of large scale deployments.  Projects
include creating new theoretical foundations for ethics in practice
\citep{GuptaEtAl2023-WhichLpNorm}, including fairness into supply
chains and energy networks \citep{HettleEtAl2021-FairReliableReconnections}, technological
progress for high-impact policy changes
\citep{GillaniEtAl2023-RedrawingAttendanceBoundaries}, and the
incorporation of IEEE Well-being Metrics into the design and
deployment of AI and optimization research.

\vspace{-0.1cm}
\section{Workforce Development}
\label{section:education}

\begin{figure}[!t]
	\centering
	\includegraphics[width=\linewidth, trim=0 0 0 5cm, clip]{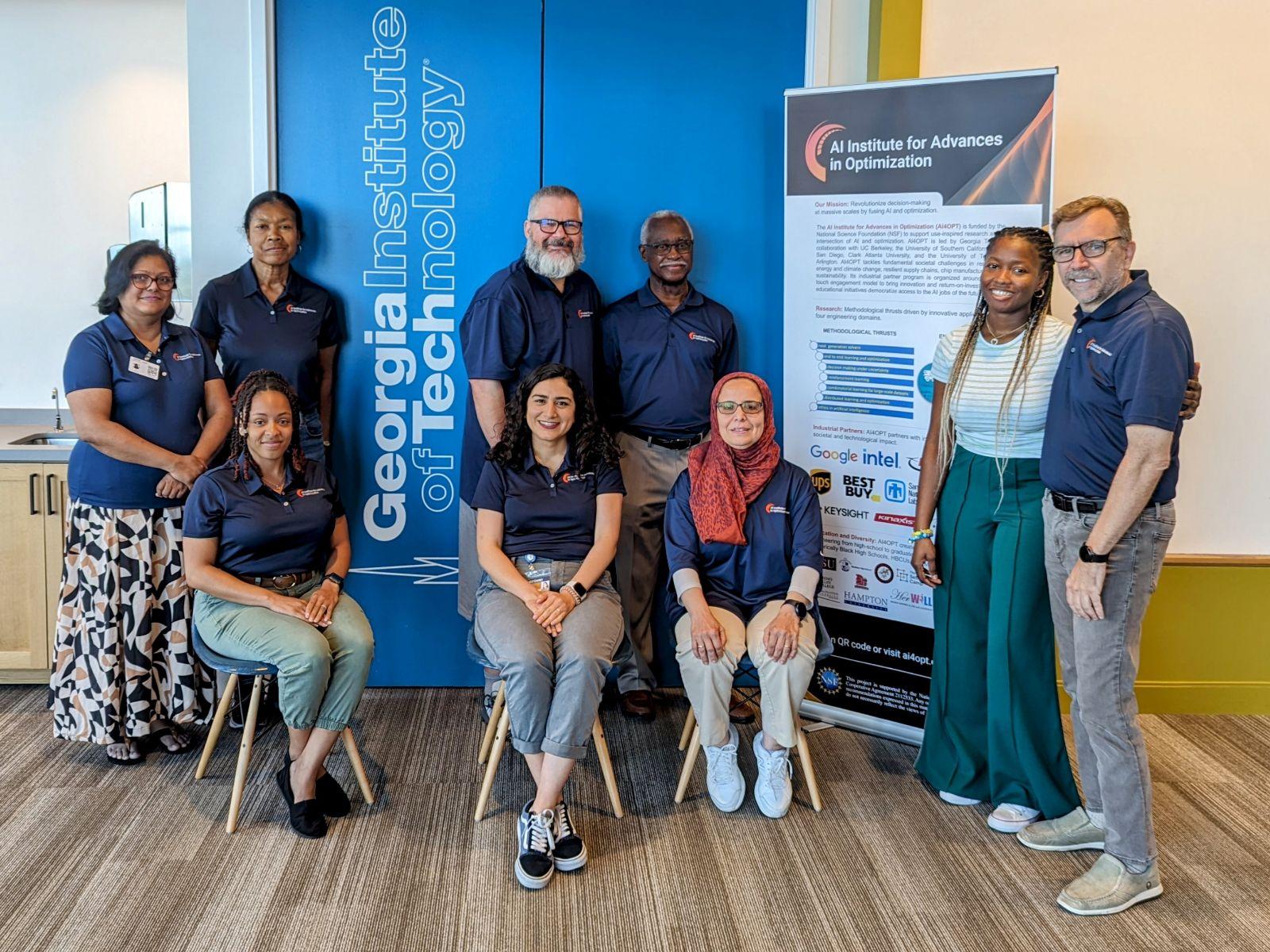}
	\caption{Faculty Training Program Cohort 1.}
	\label{fig:ftp}
\end{figure} 

The education and workforce development initiatives of \aiopt{} are
presented in detail by
\cite{PierreEtAl2022-LongitudinalEducationPrograms} and are only briefly mentioned
here. Perhaps the most distinctive feature of these programs is the
{\em``teaching the teachers'' philosophy} that permeates the
initiatives.
\aiopt{} is reaching middle and high-school students through the \href{https://www.ai4opt.org/seth-bonder-camp}{Seth Bonder summer camps} that are delivered, not only to students, but also to high-school teachers.
The camps are organized at Georgia Tech, at UC Berkeley, and online (in collaboration with Kids Teach Tech), and attract a diverse range of participants, with a large proportion of minority students and young women.
Students who succesfully complete the camp are invited to become mentors in the following year.

The \aiopt{} \href{https://www.ai4opt.org/undergraduate-education}{Faculty Training Program} (FTP) provides faculty members from HBCUs and MSIs with courses in AI, data science, and course design to create minors and majors in AI at their own institutions.
This three-year program includes a yearly 3-4 week visit to Georgia Tech and online courses throughout the year.
The program started in June 2022, and multiple FTP participants are already creating minors and majors, and are working with \aiopt{} to expand their AI offerings.

\vspace{-0.1cm}
\section{AI4OPT as a Nexus}
\label{section:nexus}

\aiopt{} acts as a nexus for AI and optimization with applications in
supply chains, energy systems, manufacturing, and sustainability. The Institute nexus is organized
around its research and educational programs, its Industrial Partner Program (IPP), its national and international collaborations, and its outreach activities. In outreach, \aiopt{} is pursuing additional partnerships, e.g., with Jackson State University in Mississipi and Huston-Tillotson University in Austin, Texas. The IPP continues to grow and covers the entire spectrum of activities in end-to-end supply chains and the planning and operations of electrical power systems. \aiopt{} features novel longitudinal internships that are piloted by the Institute at Georgia Tech, and strong collaborations with DOE national laboratories, Independent Systems Operators (ISO), and peer international institutions. A particularly exciting development is the collaboration with the AI Institute ICICLE around decentralized collaborative multimodal food supply chains with a focus on tribal communities.

%\vspace{-0.2cm}

\renewcommand*{\bibfont}{\scriptsize}
\setlength{\bibsep}{0.6pt plus 0.3ex}
\bibliographystyle{plainnat}
\bibliography{publications}

\end{document}